\documentclass{amsart}
\usepackage{amssymb,amscd}
\usepackage{times}

\usepackage[11pt,curve,matrix,arrow,frame,graph]{xy}
\newcounter{intro}

\newtheorem{theo}[intro]{Theorem}

\newtheorem{conjecture}[intro]{Conjecture}
\newtheorem{thm}{Theorem}[section]
\newtheorem{lem}[thm]{Lemma}
\newtheorem{prop}[thm]{Proposition}

\theoremstyle{remark}
\newtheorem{rem}[thm]{Remark}
\newtheorem*{merci}{Acknowledgements}

\numberwithin{equation}{section}   

\newcounter{counteroman}
\newenvironment{enumeroman}{\begin{list}{\roman{counteroman})}{\usecounter{counteroman}}}{\end{list}}

\newcommand{\cref}[1]{Corollary~\ref{#1}}

\newcommand{\R}{\mathbb{R}}
\newcommand{\C}{\mathbb{C}}
\newcommand{\B}{\mathbb{B}}
\renewcommand{\H}{\mathbb{H}}

\newcommand{\PP}{\mathbb{P}}
\newcommand{\N}{\mathbb{N}}
\newcommand{\M}{\mathbb{M}}

\newcommand{\bS}{\mathbb{S}}
\newcommand{\LL}{\mathbb{L}}

\newcommand{\cH}{\mathcal{H}}\newcommand{\mcC}{\mathcal{C}}\newcommand{\mcO}{\mathcal{O}}

\newcommand{\mcM}{\mathcal{M}}
\newcommand{\mcP}{\mathcal{P}}

\newcommand{\Gg}{\mathfrak{g}}

\newcommand{\Pg}{\mathfrak{p}}

\newcommand{\Qg}{\mathfrak{q}}
\newcommand{\Kg}{\mathfrak{k}}
\newcommand{\Ku}{\mathfrak{u}}
\newcommand{\bpm}{\begin{pmatrix}}
\newcommand{\epm}{\end{pmatrix}}

\DeclareMathOperator{\eucl}{eucl}
\DeclareMathOperator{\Ad}{Ad}

\DeclareMathOperator{\Id}{Id}
\DeclareMathOperator{\im}{Im}
\DeclareMathOperator{\Fix}{Fix}

\DeclareMathOperator{\ad}{ad}

\DeclareMathOperator{\tr}{tr}
\DeclareMathOperator{\vol}{vol}
\def\Hi#1{ {\rm Hilb}^{#1}_0(\C^2) }

\begin{document}

\title[On the QALE geometry of Nakajima's metric]
{On the QALE geometry of Nakajima's metric }

\author{Gilles Carron}
\address{Laboratoire de Math\'ematiques Jean Leray (UMR 6629), Universit\'e de Nantes, 
2, rue de la Houssini\`ere, B.P.~92208, 44322 Nantes Cedex~3, France}
\email{Gilles.Carron@math.univ-nantes.fr}

\begin{abstract} We show that on Hilbert scheme of $n$ points on $\C^2$, 
the hyperk\"ahler metric construsted by H. Nakajima via hyperk\"ahler 
reduction is the Quasi-Asymptotically Locally Euclidean (QALE in short)
 metric constructed by D. Joyce.

\end{abstract}

\maketitle

\section{Introduction}
The Hilbert scheme (or Douady scheme) of $n$ points on $\C^2$, denoted by $\Hi{n}$,
 is a crepant resolution of the quotient of
$\left(\C^2\right)^n_0=\left\lbrace q\in\left(\C^2\right)^n,\sum_j q_j=0 \right\rbrace $ 
by the action of the symmetric group $S_n$ which acts by permutation of the indices :
$$\sigma\in S_n, q\in \left(\C^2\right)^n_0,
\ \sigma.q=(q_{\sigma^{-1}(1)},q_{\sigma^{-1}(2)},...,q_{\sigma^{-1}(n)}).$$
Hence we have a map $$\pi\,:\, \Hi{n}\rightarrow \left(\C^2\right)^n_0/S_n.$$
The complex manifold $\Hi{n}$ carries a natural complex symplectic structure which 
comes from the $S_n$ invariant one of $\left(\C^2\right)^n_0$.
A compact K\"ahler manifold admitting a complex symplectic form carries in his K\"ahler
class a hyperk\"ahler metric, this is now a well-know consequence of the solution of
the Calabi conjecture by S-T. Yau (see \cite{Be}). However, $\Hi{n}$ is non compact,
for instance $\Hi{2}=T^*\PP^1(\C)$. There are many extensions of Yau's result to
non compact manifold (see for instance \cite{Bando},\cite{TY1},\cite{TY2}) and in 1999,
 D. Joyce has introduced a new class of asymptotic geometry 
 called Quasi-Asymptotically Locally Euclidean (QALE in short) ; 
 this class is the extension of the class of ALE (for Asymptotically Locally Euclidean) ;
  roughly a complete manifold $(M^d,g)$ is called ALE asymptotic to $\R^d/\Gamma$ where 
  $\Gamma\subset O(d)$ is a finite subgroup acting freely on $\bS^{d-1}$, 
  if outside a compact set $M$ is diffeomorphic to $\left(\R^d\setminus\B\right)/\Gamma$ and if
   on there the metric is asymptotic to the Euclidean metric (the precise definition requires 
   estimates between $g$ 
and the Euclidean metric). When $X^m$ is a crepant resolution of $\C^m/\Gamma$ for 
 $\Gamma\subset SU(m)$  a finite group, then roughly a K\"ahler metric on
  $X^m$ is called QALE if firstly away from the pulled back of the singular set
   the metric is asymptotic to the Euclidean one and secondly on pieces of $X^m$ 
   which (up to a finite ambiguity) are diffeomorphic to a subset of $X_A\times \Fix (A)$ 
   where $A$ is a subgroup of $\Gamma$,
and $X_A$ is a crepant resolution of $\Fix(A)^\perp/A$, then the metric is asymptotic to the sum
 of a QALE K\"ahler metric on $X_A$ and a Euclidean metric on $\Fix (A)$.  And D. Joyce has 
 proved the following (\cite{joyce_a},\cite{joyce_b}[theorem 9.3.3 and 9.3.4]):
\begin{theo} When $\Gamma\subset SU(m)$ is a finite group and $X^m\rightarrow \C^m/\Gamma$ 
is a crepant resolution , then in any K\"ahler class of QALE metric there is a unique QALE K\"ahler
 Ricci flat metric. Moreover if $\Gamma\subset Sp(m/2)$, then this metric is hyperk\"ahkler.
 \end{theo}
In particular, up to scaling, $\Hi{n}$ carries an unique hyperk\"ahler metric asymptotic 
to $\left(\C^2\right)^n_0/S_n.$

Another fruitful construction of hyperk\"ahler metric is the hyperk\"ahler quotient construction of 
N. Hitchin, A. Karlhede, U. Lindstr\"om and M. Rocek \cite{HKLR}. In fact in 1999, H. Nakajima has 
constructed a hyperk\"ahler metric on $\Hi{n}$ as a hyperk\"ahler quotient \cite{Naka}. Moreover 
H. Nakajima asked wether this metric could be recover via a resolution of the Calabi conjecture ; 
also D. Joyce said hat it is likely that QALE hyperk\"ahler metric can be explicitly constructed
 using the hyperk\"ahker quotient , but outside the case of $\Gamma\subset SU(2)=Sp(1)$
  treated by Kronheimer \cite{Kr1}, he has no examples. The main result of this paper is the following :
\begin{theo}
 On $\Hi{n}$, up a scaling, D.Joyce's and H.Nakajima's metrics coincide.
\end{theo}
It should be noted that a given complex manifold can carry two very different hyperk\"ahler metrics ;
 for instance as it has been clearly explain by C. Lebrun $\C^2$ carries two quite different K\"ahler
  Ricci flat metrics, the Euclidean one and the Taub-Nut metric which has cubic volume growth
   \cite{lebrun}.

The main evident idea of the proof of this result is to study the asymptotics of
Nakajima's metric ;
however in order to used D. Joyce unicity result,
we would need also asymptotics on the derivatives of Nakajima's metric,
this is probably possible but requires more estimates. 
Our analyze of the asymptotics of Nakajima's metric gives that Joyce and Nakajima's metrics differ 
by $O\left(\rho^{-2}\sigma^{-2}\right)$ ; where $\rho$ is the distance to a fixed point and $\sigma$ 
a regularized version of the distance to the singular set. And  in order to  used the classical argument
of S-T. Yau giving the unicity of the solution to the Calabi conjecture, we need to find a function
$\varphi$ vanishing at infinity such that, the difference between the two K\"ahler forms of
Nakajima and Joyce's metric is $i\partial\bar\partial \varphi$. D. Joyce has developed
elaborate tools to solve the equation of the type $\Delta u=f$ on QALE manifold ;
but the decay $O\left(\rho^{-2}\sigma^{-2}\right)$ is critical for this analysis.
In fact, we have circumvent this difficulty using the Li-Yau's estimates for the Green kernel
of a manifold with non negative curvature \cite{LiYau}, and we have obtained the following result,
which has a independent interest and which can be generalize to other QALE manifolds :
\begin{theo}
 Let $f$ is a locally bounded function on $\Hi{n}$, such that for some $\varepsilon>0$ satisfies
$$f=O\left(\frac{1}{\rho^{\varepsilon}\sigma^{2}}\right)$$
then the equation $\Delta u=f$ has a unique solution such that
$$u=O\left(\frac{\log(\rho+2)}{\rho^{\varepsilon}}\right).$$
\end{theo}
For further more profound results on the analysis on QALE space, there is a very interesting work of A. Degeratu and R. Mazzeo \cite{DM}.

In the physic litterature, $\Hi{n}$ is associated to the moduli space of instantons on noncommutative $\R^4$ \cite{NS}.  Our motivation for the study of the asymptotic geometry of the Nakajima's metric comes from a question of C.Vafa and E. Witten about the space of $L^2$ harmonic forms on $\Hi{n}$ endowed with the Nakajima's metric. Let $\cH^k$ be the space of $L^2$ harmonic $k-$forms on $\Hi{n}$ :
$$\cH^k=\left\lbrace \alpha\in L^2\left( \Lambda^kT^*\Hi{n}\right),\ d\alpha=d^*\alpha=0 \right\rbrace .$$
In \cite{VW}, see also the nice survey of T. Hausel \cite{Hau}, the following question is asked :
\begin{conjecture}\label{VafaWitten}
 $$\cH^k=\left\lbrace \begin{array}{ll}
\{0\}&\ {\rm if\ } k\not=2(n-1)=\dim_\R\Hi{n}\\
\im\left( H_c^k(\Hi{n})\rightarrow H^k(\Hi{n})\right) &\ {\rm if\ } k=2(n-1)
                      \end{array} \right. $$
\end{conjecture}
However, C. Vafa and E. Witten said ''Unfortunately, we do not understand the prediction of $S$-duality on non-compact manifolds precisely enough to fully exploit them.''

In fact, N. Hitchin has shown that the vanishing of the space of $L^2$ harmonics $k-$forms outside middle
 degree is a general fact for hyperk\"ahler reduction of the flat quaternionic space $\H^m$ by a compact
  subgroup of $Sp(m)$ \cite{H} ; he obtained this result with a generalization of an idea of M. Gromov
   (\cite{Gro} see also related works by J.Jost, K. Zuo and J. Mc Neal \cite{JZ},\cite{MN}). For the degree $k=2(n-1)$, the cohomology of $\Hi{n}$ is well known :
\begin{equation*}\begin{split} H_c^{2(n-1)}(\Hi{n})&\simeq\im\left( H_c^{2(n-1)}(\Hi{n})\rightarrow H^{2(n-1)}(\Hi{n})\right)\\
&\simeq  H^{2(n-1)}(\Hi{n})\simeq\R
\end{split}\end{equation*}
and a  dual class to the generator is $\pi^{-1}\{0\}$. Moreover a general result of M. Anderson says that the image of the cohomology with compact support
in the cohomology always injects inside the space of $L^2$ harmonics forms \cite{An}. Hence for the Hilbert scheme of $n$ points in $\C^2$ endowed with  Nakajima's metric we always have 
$$\dim \cH^{2(n-1)}\ge 1$$ and the conjecture \ref{VafaWitten} predicts the equality
$\dim \cH^{2(n-1)}= 1$.

There are many results on the topological interpretation of the space of $L^2$ harmonic forms on non compact manifolds but all of them requires a little on the knowledge of the asymptotic geometry (see  \cite{HHM} for results related to some prediction from string theory and \cite{L2} for a list of such results) ; the rough idea is that this asymptotic geometry would provide a certain behavior of $L^2$ harmonic forms (decay, polyhomogeneity in a good compactification) and that would imply a topological interpretation of this space with a cohomology of a compactification. With our paper \cite{L2QALE}, our main result implies :
\begin{theo}
 The Vafa-Witten conjecture \ref{VafaWitten} conjecture is true when $n=3$.
\end{theo}

The case $n=2$ can be treated by explicit computation  (see \cite{H} for clever computations).

As the Vafa-Witten conjecture is in fact more general and concerns
 the quivers varieties constructed by H.Nakajima \cite{N2}, a natural perspective is to understand
  the asymptotic geometry of the quivers varieties and  the class of Quasi-asymptotically Conical manifolds introduced by R. Mazzeo 
should be usefull \cite{M}. In a different direction it would be good to develop appropriate
 QALE tolls to settle the status of the Vafa-Witten conjecture.
\begin{merci}
It is a pleasure to thank A. Degeratu, P. Romon, R. Mazzeo, M. Singer, C. Sorger and Y. Rollin 
for interessing discussion related to this work ; a special thank is due to O. Biquard who suggested 
 that I could used the classical proof of the unicity of the solution of the Calabi conjecture in place 
 of difficult derivative estimate. This paper was finished during a stay at the MSRI,
  and I was partially supported by a joint NSF-CNRS project and the project ANR project GeomEinstein 06-BLAN-0154.
\end{merci}
\section{Nakajima's metric} In \cite{Naka}, H. Nakajima has shown that the Hilbert scheme of $n$ points in
$\C^2$ carries a natural hyperk\"ahler metric ; this metric is obtained from the Hyperk\"ahlerian quotient
construction of N. Hitchin, A. Karlhede, U. Lindstr\"om and M. Rocek \cite{HKLR}~: 
the complex vector space $$\M_n=\M:=\{(A,B,x,y)\in \mcM_n(\C)\oplus \mcM_n(\C)\oplus \C^n\oplus (\C^n)^*, \tr
A=\tr B=0\}$$
has a complex structure
$$J(A,B,x,y)=(B^*,-A^*,y^*,-x^*)$$ if we let $K=iJ$ then
$(\M,I=i,J,K:=iJ)$ becomes a quarternionic vector space ; moreover the unitary group $U(n)$ acts linearly on
$\M$ : 
if $g\in U(n)$ and $z=(A,B,x,y)\in \M$ then
$$g.z=\left(gAg^{-1},gBg^{-1},gx,yg^{-1}\right).$$
The real moment map associated to this action is 
$$\mu(A,B,x,y)=\frac{1}{2i}\left(\left[A,A^*\right]+\left[B,B^*\right]+xx^*-yy^*\right)\in \Ku(n).$$
If $h\in \Ku(n)$ and $z=(A,B,x,y)\in \M$ we let
$$l_z(h)=\left.\frac{d}{dt}\right|_{t=0}e^{th}.z=\left([h,A], [h,B], hx,-yh\right).$$
By definition, we have for
$z\in \M, \delta z\in T_z\M\simeq \M$ :
$$\langle d\mu(z)(\delta z),h\rangle=\langle i l_z(h),\delta z\rangle$$
The action of $GL_n(\C) $ on $\M$ preserves the
complex symplectic form
$$\omega_\C(z, z')=\tr \left( A .B'- B .A' \right)
+y'(x)-y(x'),$$ 
and the associated complex moment map is :
$$\mu_\C(A,B,x,y)=[A,B]+xy\in \mcM_n(\C).$$
Let $t>0$ and defined
$$\LL_t(n)=\LL_t:=\mu^{-1}\left\{\frac{t}{2i}\right\}\cap \mu_\C^{-1}\{ 0\}, $$
then the map $${\mu\!\!\!\mu}:=(\mu,\mu_\C)\,:\, \M\rightarrow \Ku(n)\oplus \mcM_n(\C)$$ is a submersion near
$\LL_t$ and $U(n)$ acts freely on it, hence the quotient $\H_t:=\LL_t/U(n)$ is a smooth manifold, 
this manifold is endowed with the Riemannian metric $g_N$ which makes the submersion
$\LL_t \rightarrow \LL_t/U(n)$ Riemannian. By definition the tangent space of $U(n)z$ is naturally isometric
to the orthogonal of the space
$$\im l_z\oplus I\im l_z\oplus J \im l_z\oplus K\im l_z.$$
In particular, $\H_t$ is endowed with a quaternionic structure which is in fact integrable ; hence
the metric $g_n$ is hyperk\"ahler hence K\"ahler and Ricci flat.
\subsection{Some remarks} 
Because for $\lambda>0$, we have $\lambda\LL_t=\LL_{\lambda^2 t}$, all the spaces
$\{\H_t\}_{t>0}$ are isomorphic and their Riemannian metrics are proportional.

 For $t=0$, the quotient $\LL_0/U(n)$ is not a smooth manifold. It is easy to show
that
$$(A,B,x,y)\in \LL_{0}  \Leftrightarrow \left(x=0,y=0, [A,B]=[A,A^*]=[B,B^*]=0\right) ;$$
hence to $(A,B,0,0)\in \LL_{0}$ we can associated their joint spectrum
$$S_n.\left((\lambda_1,\mu_1),(\lambda_2,\mu_2),...,(\lambda_n,\mu_n)\right) \in \left(\C^2\right)_0^n/S_n$$
where $\left(\C^2\right)_0^n:=\{(q_1,...,q_n\in\left(\C^2\right)^n, \sum_j q_j=0\}$ and the symmetric group $S_n$ acts on 
$\left(\C^2\right)_0^n$ by permutation of the indices.
We get an isomorphism (in fact an isometry)
$$\LL_{0}/U(n)\simeq \left(\C^2\right)_0^n/S_n.$$
In fact for $t>0$, we still have
\begin{equation}\label{y0}
(A,B,x,y)\in \LL_{t}  \Rightarrow y=0.\end{equation} 
Hence for $z:=(A,B,x,0)\in \LL_t$, the joint spectrum
of $(A,B)$ is still defined and we can defined
$$\pi(U(n)z)=S_n.\left((\lambda_1,\mu_1),(\lambda_2,\mu_2),...,(\lambda_n,\mu_n)\right)  \in \left(\C^2\right)_0^n/S_n$$
where $(\lambda_1,\mu_1),(\lambda_2,\mu_2),...,(\lambda_n,\mu_n)$ are such that for a $g\in U(n)$,
the matrix $gAg^{-1}$ (resp. $gBg^{-1}$) is upper triangular with diagonal
$(\lambda_1,\lambda_2,..,\lambda_n)$ (resp. $(\mu_1,\mu_2,..,\mu_n)$ ).
$\H_t $ is  isomorphic to the Hilbert scheme of $n$ points\footnote{with the
center of mass removed.} in $\C^2$ $\Hi{n}$.
The map $\pi \,:\,\Hi{n}  \rightarrow \left(\C^2\right)^n_0/S_n$ is in fact a crepant resolution of $\left(\C^2\right)_0^n/S_n$.
\begin{rem}\label{yzero}
We also remark that if $v=(\delta A, \delta B,\delta x,0)\in T_\zeta\LL_t$
 is orthogonal to the range of $l_\zeta$, then $Jv$ is also in $T_\zeta\LL_t$, and hence $\delta x=0$.
\end{rem}
\subsection{The geometry of $\Hi{2}$.}
As an example, we look at the geometry of $\Hi{2}$. Let
$z=(A,B,x,0)\in \mcM_2(\C)\oplus\mcM_2(\C)\oplus \C^2\oplus \left(\C^2\right)^*$ such that
$\tr A=\tr B=0$ and
$$\left\{\begin{array}{l}
\left[A,A^*\right]+\left[B,B^*\right]+xx^*=t \Id\\
\left[A,B\right]=0\end{array}
\right. .$$

When $\det A\not=0$ or $\det B\not=0$ then we can find a $g\in U(2)$ such that
$$gAg^{-1}=\bpm\lambda&a\\0&-\lambda\epm ,gBg^{-1}=\bpm\mu&b\\0&-\mu\epm $$
then let $g(x)=(x_1,x_2)$.
The equation $\left[A,B\right]=0$ implies that there is a number $\rho$ such that
$a=\lambda \rho $ and $b=\mu   \rho$. Then the remaining equations are
for $R^2:=|\lambda|^2+|\mu|^2$
$$\left\{\begin{array}{l}
|\rho|^2R^2+|x_1|^2=t\\
-|\rho|^2R^2+|x_2|^2=t\\
-2R^2\rho+x_1\overline{x_2}=0\\
\end{array}
\right.$$
We can always choose $g$ such that $x_1,x_2\in \R_+$
then we obtain
$\rho^2=\sqrt{4+\frac{t^2}{R^4}}-2$
Hence
$$\rho=\frac{t}{2R^2}+O\left(\frac{1}{R^6}\right)$$
if $x_1=\sqrt{2t}\sin(\phi),x_2=\sqrt{2t}\cos(\phi)$, then
$t\cos(2\phi)=\rho^2R^2$ and $t\sin(2\phi)=2\rho R^2$
Hence $\phi=\frac{\pi}{4}+O\left(\frac{1}{R}\right)$ and
$x_1=\sqrt{\frac{t}{2}}+O\left(\frac{1}{R}\right)$, $x_2=\sqrt{\frac{t}{2}}+O\left(\frac{1}{R}\right)$.
Hence for $(\lambda,\mu)\in \C^2\setminus\{0\}\simeq \left(\C^2\right)^2_0$ we have found
$$z(\lambda,\mu)= 
\left(\bpm\lambda&\lambda\rho(R)\\0&-\lambda\epm ,\bpm\mu&\mu\rho(R)\\0&-\mu\epm,\,
x(R),0\,\right)\in \LL_t.$$
Moreover, $z(\lambda,\mu)$ and $z(\lambda',\mu')$ are in the same $U(2)$ orbit if and only if
$(\lambda,\mu)=\pm(\lambda',\mu')$ ; hence we have a map
$$z\,:\, \left(\C^2\setminus\{0\}\right)/\{\pm Id\}\rightarrow \LL_t/U(2).$$
From the exact value of $z$, we can show that
$$z^*g_N=2 \left[ |d\lambda|^2+|d\mu|^2\right]+O\left(\frac{1}{R^4}\right).$$
This shows that $(\Hi{2},g_N)$ is a hyperk\"ahler metric which is Asymptotically
Locally Euclidean asymptotic to $\C^2/\{\pm Id\}$. These manifolds has been classified by 
Kronheimer \cite{Kr2}, so that in this case Nakajima's metric is the Eguchi-Hansen metric on
$T^*\PP ^1(\C)$.
\subsection{A last useful remark} A priori, it is not clear wether the above map $z$ is holomorphic, this is
in fact true as a consequence of the following useful lemma :
\begin{lem}\label{holo}
Suppose that a compact Lie group  $G$ acts on $\H^m$ by quaternionic linear maps and let
$\mu\!\!\!\!\mu\,:\, \H^m\rightarrow \Gg^*\otimes \im\H$ be the associated moment map. Assume that 
for some $\zeta=(\zeta_\R,\zeta_\C)\in\Gg^*\otimes \im\H$ the hyperk\"ahler quotient $Q:=\mu\!\!\!\!\mu^{-1}\{\zeta\}/G$ is well
defined. When $X$ is a complex manifold and 
 $\Psi\,:\, X\rightarrow \mu\!\!\!\!\mu^{-1}\{\zeta\}$ is a smooth map such that
locally $$\Psi(x)=g(x)\tilde\Psi(x)$$
where $g\,:\, X\rightarrow G^{\C}$ is smooth and $\tilde \Psi\,:\, X\rightarrow \mu_\C^{-1}\{\zeta_\C\}$
is holomorphic, then the induced map $\bar\Psi\,:\, X\rightarrow Q$ is also holomorphic.
\end{lem}
\proof
If $q\in \H^m$ let $P_q$ be the orthogonal projection onto the orthogonal of
$\im l_q\oplus I\im l_q\oplus J\im l_q\oplus K\im l_q=\im l^\C_q\oplus J\im l^\C_q$
where $l_q\,:\, \Kg\rightarrow \H^q$ is defined as before by
$$l_q(h)=\left.\frac{d}{dt}\right|_{t=0}e^{th}.q=h.q\,.$$
We must show that if $x\in X$, then for $q:=\Psi(x)$ :
$$P_q\big(d\Psi(x)(Iv)\big)=IP_q\big(d\Psi(x)(v)\big).$$
But  $\dot g(x)=dg(x)(Iv)\in G^\C$ we have
$$d\Psi(x)(Iv)=\dot g(x).q+g(x).d\tilde\Psi(x)(Iv)=l^\C_q(\dot g(x))+g(x).d\tilde\Psi(x)(Iv)$$
By definition  $P_q(l^\C_q(\dot g(x)))=0$ and because $g(x)$ and $P_q$ are complex linear :
$$P_q\big(d\Psi(x)(Iv)\big)=P_q\big(g(x).Id\tilde\Psi(x)(v)\big)=IP_q\big(d\Psi(x)(v)\big).$$
\endproof 
\section{Joyce's metric}
In \cite{joyce_a,joyce_b}, D. Joyce has build many new K\"ahler, Ricci flat metrics on some crepant
 resolution of quotient of
$\C^m$ by a finite subgroup of $SU(m)$ ; his construction relies upon the resolution of a Calabi-Yau
 problem
for a certain class of asymptotic geometry which is called QALE for
 Quasi Asymptotically Locally Euclidean. 
We will follow the presentation of D. Joyce for the Hilbert Scheme of $n$ points
 on $\C^2$ and we will then
describe the asymptotic geometry of these QALE metrics on $\Hi{n}$.
\subsection{The local product resolution of $\Hi{n}$.}
If $\Pg=(I_1,I_2,...,I_k)$ is a partition of $\{1,2,..,n\}$
\footnote{the $I_l$'s are disjoint and their reunion
is $\{1,2,...,n\}$.}, the $I_l$'s are called the cluster of $\Pg$. We will denote
$$V_\Pg=\{q\in \left(\C^2\right)^n_0,\forall l\in\{1,...,k\}, \forall i,j\in I_l : q_i=q_j\}$$
and
$A_\Pg=\{\gamma\in S_n,\ \gamma q=q\,\forall q\in V_\Pg\}\simeq S_{n_1}\times S_{n_2}\times...\times S_{n_k}$
where $n_l=\#I_l.$ Then
$$W_\Pg=V_\Pg^\perp\simeq \bigoplus_{l=1}^k \left(\C^2\right)^{n_l}_0. $$
Let $m_\Pg={\rm codim}_\C V_\Pg=\dim_\C W_\Pg=2(n-l(\Pg))$ where $l(\Pg)=k$. 
The set $\mcP_n$ of partitions of $\{1,2,...,n\}$ has the following partial order :
$$\Pg\leq \Qg \Leftrightarrow V_{\Qg}\subset V_\Pg \Leftrightarrow W_\Pg\subset W_\Qg$$
Hence $\Pg\leq \Qg$ if and only if $\Pg$ is a refinement of $\Qg$ : i.e. if
$\Qg=(J_1,J_2,...,J_r)$, then there are partitions $(I_{l,1},I_{l,2},...,I_{l,n_l})$ of
$J_l=I_{l,1}\cup...\cup I_{l,n_l}$ such that the cluster of $\Pg$ are the
$ I_{l,j}$'s.
The smallest partition is $\Pg_0=\{1\}\cup\{2\}\cup...\cup\{n\}$ with $V_{\Pg_0}=\left(\C^2\right)^n_0$,
 the largest
partition is $\Pg_\infty=\{1,2,...,n\}$ with $V_{\Pg_\infty}=\{0\}.$
The fundamental partitions are the $\Pg_{i,j}$ such that
$$\Pg_{i,j}=\big(\{i,j\},\{k_1\},\{k_2\},...,\{k_{n-2}\}\big)$$
with $\{1,2,...,n\}\setminus\{i,j\}=\{k_1,k_2,...,k_{n-2}\},$
then $V_{i,j}:=V_{\Pg_{i,j}}=\{q\in \left(\C^2\right)^n_0, q_i=q_j\}$.
We have for any partition $\Pg\not=\Pg_0$
$$V_\Pg=\cap_{\Pg_{i,j}\leq \Pg}V_{i,j}$$
We will also denote
$\Delta_\Pg=\{(i,j)\in \{1,2,...,n\}^2,\ \Pg_{i,j}\not\leq \Pg\}$ and
$\Delta^c_\Pg=\{(i,j)\in \{1,2,...,n\}^2,\ \Pg_{i,j}\leq \Pg\}$.
The singular locus of $\left(\C^2\right)^n_0/S_n$ is the quotient of the generalized diagonal
$$S=\left(\bigcup_{\Pg\not=\Pg_0}V_p\right)/S_n=\left(\bigcup_{i,j} V_{i,j}\right)/S_n.$$
Finally let 
$$S_\Pg=\left(\bigcup_{(i,j)\in \Delta_\Pg}V_{i,j}\right)/A_\Pg$$ and for $R>0$, let $T_\Pg$ be the
$R$-neighborhood of $S_\Pg$ is $\left(\C^2\right)^n_0/A_\Pg$ :
$$T_\Pg:=\{q\in \left(\C^2\right)^n_0,\  \exists (i,j)\in \Delta_\Pg\  |q_i-q_j|<R\}/A_\Pg.$$
The resolution $\pi\,:\,\Hi{n}\rightarrow \left(\C^2\right)^n_0/S_n$  is a local product resolution ;
indeed there is a resolution of $W_\Pg/A_\Pg$ namely
$$\pi_\Pg\,:\,\Hi{\Pg}:=\prod_{l=1}^k\Hi{n_l}\rightarrow W_\Pg/A_\Pg$$ such that for
$U_\Pg=(\pi_\Pg\times \Id)^{-1}(T_\Pg)\subset\Hi{\Pg}\times V_{\Pg} $ and
$\phi_\Pg\,\:\left(\C^2\right)^n_0/A_\Pg\rightarrow \left(\C^2\right)^n_0/S_n$ the natural map
 there is a local biholomorphism onto his image $\psi_\Pg\,:\,U_\Pg\rightarrow
\Hi{n}$ for which the following diagram is commutative :

$$\xymatrix{
{\Hi{\Pg}\times V_{\Pg}\setminus U_{\Pg}}\ar[d]^{\pi_\Pg\times\Id}\ar[r]^{\ \ \psi_\Pg} &{\Hi{n}}\ar[d]^{\pi}\\
{{\left(\C^2\right)^n_0/A_\Pg}\setminus T_\Pg}\ar[r]^{\ \ \phi_\Pg}&{\left(\C^2\right)^n_0/S_n}
} $$
In the hyperk\"ahlerian quotient description, the local biholomorphism $\psi_\Pg$ is given as follows
identifying $V_\Pg$ with $\left(\C^2\right)^k_0$, if we
let
 $$\zeta=\left((A_1,B_1,x_1,0), (A_2,B_2,x_2,0),...,(A_k,B_k,x_k,0)\right) \in \prod_{j=1}^k \LL_t(n_j)$$
and
$\eta=\left((\lambda_1,\mu_1),(\lambda_2,\mu_2),...,(\lambda_k,\mu_k)\right)\in \left(\C^2\right)^k_0\setminus U_\Pg$,
 we associated to $(\zeta,\eta)$, the vector $(A,B,x,0)\in \M(n)$
  such that $A$ and $B$ are block diagonal  with respective diagonal 
  $(A_1+\lambda_1,A_2+\lambda_2,...,A_k+\lambda_k)$ and $(B_1+\mu_1,B_2+\mu_2,...,B_k+\mu_k)$ 
  and $x=(x_1,x_2,...,x_k)$ then 
$\psi_\Pg\left( \left(U(n_1)\times U(n_2)\times...\times U(n_k)\right).\zeta,\eta\right)$ 
is the set of points, in the $GL_n(\C)$-orbit of  $(A,B,x,0)$, satisfying the real moment map equation
 (see the part 4 for more details).
\subsection{QALE metric on $\Hi{n}$.}
We introduce several functions of distance's type on 
$\left(\Hi{\Pg}\times V_{\Pg}\right)\setminus U_{\Pg}$, if $z\in \Hi{\Pg}\times V_{\Pg}\setminus U_{\Pg}$ and
$v=(\pi_\Pg\times\Id)(z)$, we note
$$\mu_{\Pg,\Qg}(z)=\inf_{\gamma\in A_\Pg} d(\gamma.v,V_\Qg)=d\left(v,(A_\Pg V_\Qg)/A_\Pg\right) $$ and
$$\nu_\Pg(z)=1+\inf_{p\not=p_0}\mu_{\Pg,\Qg}(z)$$
Then a Riemannian metric $g$ on $\Hi{n}$ is called QALE (asymptotic to $\left(\C^2\right)^n_0/S_n$) if for each partition
$\Pg$ there is a metric $g_\Pg$ on $\Hi{\Pg}$  such that
for all $l\in \N$ :
\begin{equation}\label{QALE1}
\nabla^l\left(\psi_\Pg^*g-(g_\Pg+\eucl_{V_\Pg})\right)=\sum_{\Qg\not\leq \Pg}
O\left(\frac{1}{\nu_\Pg^{2+l}\mu_{\Pg,\Qg}^{2m_\Qg-2}}\right)\end{equation}
However, if $\Qg\not\leq \Pg$ there is always a $(i,j)\in \Delta_\Pg$ such that $\Pg_{i,j}\not\leq \Pg$ 
and $\Pg_{i,j}\leq \Qg$ , therefore  $$\mu_{\Pg,\Pg_{i,j}}^{2m_{\Pg_{i,j}}-2}=\mu_{\Pg,\Pg_{i,j}}^{2}\leq
\mu_{\Pg,\Qg}^{2m_\Qg-2}$$
If we introduce
$\rho_\Pg(z)=\inf_{(i,j)\in \Delta_\Pg} \mu_{\Pg,\Pg_{i,j}}$ then in fact
for $v=(\pi_\Pg\times \Id)(z)\in \left(\C^2\right)^n_0/A_\Pg$ we have
$$\rho_\Pg(z)=\inf_{(i,j)\in \Delta_\Pg} \left|v_i-v_j\right|$$
The asymptotic \ref{QALE1} are equivalent to
\begin{equation}\label{QALE2}
\nabla^l\left(\psi_\Pg^*g-(g_\Pg+\eucl_{V_\Pg})\right)=
O\left(\frac{1}{\nu_\Pg^{2+l}\rho_\Pg^2}\right)\end{equation}

We can introduce two functions of distance's type :
when $z\in \Hi{n}$ and $\pi(z)=(v_1,v_2,...v_n)\in \left(\C^2\right)^n_0/S_n$, then we let
$$\rho(z)=\sqrt{\sum_{i<j} |v_i-v_j|^2}, $$
and
$$\sigma(z)=\inf_{i\not=j}\{ |v_i-v_j|\}+1$$
If $\Pg$ is a partition of $\{1,2,..,n\}$, and $\epsilon,\tau,R$ are positive real numbers, then we introduce
\begin{equation*}\begin{split}
\check\mcC_{\Pg_0}=\{(v_1,...,v_n)\in \left(\C^2\right)^n_0/S_n, &{\rm\  such\ that\ } |v|>R\\
&{\rm\  and\ }\forall i\not=j,\ 
|v_i-v_j|>\varepsilon |v|\}\end{split}\end{equation*}
\begin{equation*}\begin{split}
\check\mcC_\Pg=\{(v_1,...,v_n)\in \left(\C^2\right)^n_0/A_\Pg, &{\rm\  such\ that\ } |v|>R,\\
&\forall (i,j)\in \Delta_\Pg\ |v_i-v_j|>\sqrt{\frac{2}{n(n-1)}}\,|v|\\
&{\ \rm\  and\ } \forall (i,j)\in \Delta^c_\Pg,\ |v_i-v_j|<2\epsilon |v|\}\end{split}\end{equation*}
It is clear that if $\epsilon$ is small enough  then the
$$\left(\C^2\setminus R\B\right)^n_0/S_n=\cup_{\Pg}\phi_\Pg(\check\mcC_\Pg).$$
Moreover on  $\mcC_\Pg:=(\pi_\Pg\times\Id)^{-1}\left(\check\mcC_\Pg\right)$, the asymptotic $\ref{QALE2}$ are
\begin{equation}\label{QALE3}
\nabla^l\left(\psi_\Pg^*g-(g_\Pg+\eucl_{V_\Pg})\right)=
O\left(\frac{1}{\sigma^{2+l}\rho^2}\right)\end{equation}
\begin{rem} It can be shown that if all metric $g_{\Pg}$ are QALE and if the estimate \ref{QALE3} is satisfied
then $g$ is also QALE.\end{rem}
\subsection{Joyce's result.}
The result of D. Joyce concerning the Hilbert scheme of $n$ points on $\C^2$ is the following :
\begin{thm}Up to scaling, $\Hi{n}$ has a unique hyperk\"ahler metric which is QALE asymptotic to
$\left(\C^2\right)^n_0/S_n$.
\end{thm}
\section{Asymptotic of Nakajima's metric}\subsection{Induction's hypothesis}
In this part, we will prove the following result by induction on $n$ :
\begin{enumeroman}
\item On $\Hi{n}$, Nakajima's metric $g_N$ satisfies the estimate (\ref{QALE3}) for $l=0$, more precisely if
$g_\Pg$ is the sum  of Nakajima's metric on $\Hi{\Pg}$, then for all partition $\Pg$ then for $\epsilon>0$
small enough and $R$ large enough, we have on $(\pi_\Pg\times\Id)^{-1}\left(\check\mcC_\Pg\right)$
$$\psi_\Pg^*(g_N)-g_\Pg+\eucl_{V_\Pg}=O\left(\frac{1}{\sigma^2\rho^2}\right)$$
\item There is a constant $C$ such that if $z=(A,B,x,0)\in \LL_t$ then
$$\forall h\in \Ku_n, \ \|l_z(h)\|^2=\|[h,A]\|^2+\|[h,B]\|^2+\|hx\|^2\geq C\|h\|^2$$
\item There is a constant $M$ such that for all  $z\in \LL_t$ and $(\delta A, \delta B, 0,0)\in
T_z\LL_t$ orthogonal to $\im l_z$ then
$$\left\|\left[\delta A,\delta A^*\right]\right\|+\left\|\left[\delta B,\delta B^*\right]\right\|\le
\frac{M}{\sigma^2}\left(\|\delta A\|^2+\|\delta B\|^2\right).$$
\end{enumeroman}

It is easy to check these three conditions for $\Hi{2}$ thanks to the explicit description of
$\LL_t$ in this case. So we now assume that these induction hypothesis are true for all $m<n$.
\subsection{The case of well separated points} We first examine the easiest case corresponding 
to $\Pg_0$. More precisely, we consider $\textbf{q}=(q_1,q_2,...,q_n)\in \left(\C^2\right)^n_0$ such that
for all $i\not=j$, then $|q_i-q_j|>R$ ($R$ will be chosen large enough), the set of such $q$'s will be
denote by $\mcO_0$.
If $q_j=(\lambda_j,\mu_j)$ we search a solution $z=(A,B,x,0)\in\M$ of the equation
\begin{equation}
\label{moments1}
\left\{\begin{array}{l}
  \left[A,A^*\right]+\left[B,B^*\right]+xx^*=t\Id\\
  \left[A,B\right]=0 \\ \end{array}\right.
\end{equation}
Where $A$, $B$ are upper triangular matrices with respective diagonals
$(\lambda_1,\lambda_2,...,\lambda_n)$ and $(\mu_1,\mu_2,...,\mu_n)$ and upper diagonal coefficients
$\textbf{a}=(a_{i,j}),\textbf{b}=(b_{i,j}) $.
We obtain the following equation\footnote{with the convention that $a_{i,j}=b_{i,j}=0$ if $j\le i$.}
 for the $(i,j)$ coefficients of the equation (\ref{moments1}) :
\begin{equation}
\label{moment2}\left\{\begin{array}{l}
(\bar\lambda_i-\bar\lambda_j)a_{i,j}+(\bar\mu_i-\bar\mu_j)b_{i,j}+
\sum_k\left[ \bar a_{k,i}a_{k,j}+\bar b_{k,i}b_{k,j}- 
a_{i,k}\bar a_{j,k}- b_{i,k}\bar b_{j,k}\right]=x_i\bar x_j\\
-(\bar\mu_i-\bar\mu_j)a_{i,j}+(\lambda_i-\lambda_j)b_{i,j}+
\sum_k\left[ a_{i,k} b_{k,j} -b_{i,k}a_{k,j}\right]=0
\end{array}\right.\end{equation}
And the equation for the diagonal coefficient $(i,i)$ of (\ref{moments1}) gives :
\begin{equation}
\label{moment3}
\sum_k \left[|a_{i,k}|^2- |a_{k,j}|^2+|b_{i,k}|^2- |b_{k,j}|^2 \right] +|x_i|^2=t
\end{equation}
We let $R_{i,j}=\sqrt{|\lambda_i-\lambda_j|^2+|\mu_i-\mu_j|^2}$ and
\begin{equation}\left\{\begin{array}{l}
x_i^0=\sqrt{t}\\
a_{i,j}^0=(\lambda_i-\lambda_j)\frac{t}{R_{i,j}^2}\\
b_{i,j}^0=(\mu_i-\mu_j)\frac{t}{R_{i,j}^2}\end{array}\right.\end{equation}
Then if we write the equations (\ref{moment2},\ref{moment3}) in the synthetic form
$$F(\textbf{q},\textbf{a},\textbf{b},\textbf{x})=0$$
where $F\,:\, \left(\C^2\right)^n_0\times \C^{n(n-1)/2}\times\C^{n(n-1)/2}\times \C^n
\rightarrow \C^{n(n-1)/2}\times\C^{n(n-1)/2}\times \C^n$. We have 
$$F\left(\textbf{q},\textbf{a}^0,\textbf{b}^0,\textbf{x}^0\right)=O\left(\sigma^{-2}\right)$$
Moreover it is easy to check that when $\sigma$ is large enough, 
the partial derivative in the last three argument
$D_{(\textbf{a},\textbf{b},\textbf{x})}
F\left(\textbf{q},\textbf{a}^0,\textbf{b}^0,\textbf{x}^0\right)$
is invertible and the norm of the inverse is uniformly bounded. The map $F$ 
being polynomial of degre $2$ in its arguments, the implicit function theorem implies that the equations
(\ref{moment2},\ref{moment3}) have a unique solution such that
\begin{equation}\label{OffDia}
(\textbf{a},\textbf{b},\textbf{x})(q)=
(\textbf{a}^0,\textbf{b}^0,\textbf{x}^0)+O\left(\sigma^{-2}\right).\end{equation}
Moreover
$D_q(\textbf{a},\textbf{b},\textbf{x})(q)=O\left(\sigma^{-2}\right)$.
We have then build a map
\begin{equation*}\begin{split}&\widehat\Psi_0\,:\, \mcO_0\rightarrow \LL_t\\ 
&\ \ \ \ q\mapsto (A(\textbf{q}), B(\textbf{q}),x(\textbf{q}),0)\\
\end{split}\end{equation*}
Moreover $\widehat \Psi_0(\textbf{q})$ and $\widehat \Psi_0(\textbf{q'})$ live in the same
$U(n)$-orbit if and only if $\textbf{q}$ and $\textbf{q'}$ live in the same $S_n$-orbit hence
$\widehat \Psi_0$ induces a map
$$\Psi_0\,:\, \mcO_0/S_n\rightarrow \Hi{n}$$ which is holomorphic according to the lemma (\ref{holo}).
With (\ref{OffDia}), we have
\begin{equation}\label{est1}
\left|d\Psi_0(\textbf{q}).v\right|^2=|v|^2+O\left( \frac{|v|^2}{\sigma^4} \right)\end{equation}
The first term comes from the diagonals of $A$ and $B$ the second one from
the off-diagonal terms and the derivative of $\textbf{x}$.
In order to check the point $i)$ of the induction hypothesis
we must show that
$$\left|d\Psi_0(\textbf{q}).v\right|^2-\left|\Pi_z\left(d\Psi_0(\textbf{q}).v\right)\right|^2=
|v|^2+O\left( \frac{|v|^2}{\sigma^4} \right)\ .$$
Where if $\widehat\Psi_0(q)=z\in \LL_t$, $\Pi_z$ is the orthogonal projection onto the space
$\im l_z$.
But by construction, if $X\in \im l_z$ then $IX$ is normal to $T_z\LL_t$ hence
$d\Psi_0(\textbf{q}).(Iv)\perp I X$ in particular
$\Pi_z\left(I.d\Psi_0(\textbf{q}).(Iv)\right)=0;$
Hence
\begin{equation}\label{est2}\Pi_z\left(d\Psi_0(\textbf{q}).v\right)
=\Pi_z\left(d\Psi_0(\textbf{q}).v+Id\Psi_0(\textbf{q}).Iv\right)=
2\Pi_z\left(\bar\partial\Psi_0(\textbf{q}).v\right).\end{equation}
But by construction
\begin{equation}\label{est3}\left|\bar\partial\Psi_0(\textbf{q}).v\right|^2
=\left|\bar\partial\textbf{a}\right|^2+
\left|\bar\partial\textbf{b}\right|^2+\left|\bar\partial\textbf{x}\right|^2
=O\left( \frac{|v|^2}{\sigma^4} \right)
\end{equation}
The assertion i) of the induction hypothesis $i)$ follows from the estimates
(\ref{est1},\ref{est2},\ref{est3}).

For the induction hypothesis $ii)$, we have for 
$z=\widehat \Psi_0(\textbf{q})$ and $h=(h_{i,j})\in \Ku_n$ :
$$\|l_z(h)\|^2\ge \frac12 \left[\sum_{i,j}R^2_{i,j} |h_{i,j}|^2+t\sum_{i} |h_{i,i}|^2\right]-C \sigma^{-2}\|h\|^2$$
Hence if $R$ is chosen large enough the induction hypothesis $ii)$ hold on $\mcO_0$.

Now we check the induction hypothesis $iii)$ let 
$$(\delta A,\delta B,0,0)=d\Psi_0(\textbf{q}).v-\Pi_z\left(d\Psi_0(\textbf{q}).v\right)$$
we have just said that
$$\|\Pi_z\left(d\Psi_0(\textbf{q}).v\right)\|=O\left(\sigma^{-2}\right) |v|.$$
Hence the off-diagonal part of $\delta A$ and $\delta B$ are bounded by $O\left(\sigma^{-2}\right) |v|$, this implies that
$$\left\|\left[\delta A,\delta A^*\right]\right\|^2+\left\|\left[\delta B,\delta B^*\right]\right\|^2\le
O\left(\sigma^{-4}\right)|v|^4.$$

\subsection{The general case}We examine now the region $\mcC_\Pg$ associated to another partition
$\Pg\not=\Pg_0$., we can always assume that 
$$\Pg=\left(\{m_0=1,...,m_1\},\{m_1+1,m_2\},...,\{m_{k-1}+1,..,n=m_k\}\right),$$
let $n_l=m_l-m_{l-1}$.
We consider the set $\mcO_\Pg$ of
$$\left(\textbf{q},\textbf{A},\textbf{B},\textbf{x}\right)\in
\left(\C^2\right)_0^k\times\bigoplus_{j=1}^k\mcM_{n_j}(\C)\times
\bigoplus_{j=1}^k\mcM_{n_j}(\C)\times\bigoplus_{j=1}^k\C^{n_j}$$
such that if $\textbf{q}=(q_1,q_2,...,q_k)$ then for all $i\not=j$ then
$|q_i-q_j|>\sqrt{\frac{1}{n(n-1)}} |\textbf{q}|$ and $|\textbf{q}|\geq R$
and if $\textbf{A}=(A_1,A_2,...,A_k),\ \textbf{B}=(B_1,B_2,...,B_k),\textbf{x}=(x_1,x_2,...,x_k)$ then each
$(A_j,B_j,x_j)$ satisfies 
$\tr A_j=\tr B_j=0$ and the moment map equation :
$$\left\{\begin{array}{l}
\left[A_j,A_j^*\right]+\left[B_j,B_j^*\right]+x_jx_j^*=t \Id_{n_j}\\
\left[A_j,B_j\right]=0\end{array}
\right.$$
and moreover 
$$\sup_j \|A_j\|^2+\|B_j\|^2\le \tau^2 |\textbf{q}|^2$$
We will search a solution $z=(A,B,x,0)$ of the moment map equation which is approximatively
\begin{equation*}
 \begin{split}
 &A\simeq\bpm A_1+\lambda_1&0&\hdots&0\\
0&A_2+\lambda_2&\ddots&\vdots\\
\vdots&\ddots&\ddots&\vdots\\
0&\hdots&\hdots&A_k+\lambda_k\epm ,\\
& B\simeq\bpm B_1+\mu_1&0&\hdots&0\\
0&B_2+\mu_2&\ddots&\vdots\\
\vdots&\ddots&\ddots&\vdots\\
0&\hdots&\hdots&B_k+\mu_k\epm,\\
& x\simeq(x_1,x_2,...,x_k)\end{split}
\end{equation*}

We first fix some $\zeta=\left(\textbf{q},\textbf{A},\textbf{B},\textbf{x}\right)\in \mcO_\Pg$
and we search a $z_0=(A^0,B^0,x^0,0)$ where if $q_j=(\lambda_j,\mu_j)$ then 
$x^0=(x_1,x_2,...,x_k)$, $A^0$ (resp. $B^0$) is upper block triangular  with diagonal
$(A_1+\lambda_1,A_2+\lambda_2\Id,...,A_k+\lambda_k)$ (resp. $(B_1+\mu_1,B_2+\mu_2,...,B_k+\mu_k)$)
and $\mu(z_0),\mu_\C(z_0)$ are block diagonal . Hence we search matrices $A_{i,j},B_{i,j}\in \mcM_{n_i,n_j}(\C), i<j$ such
that for all $i<j$ :
\begin{equation}\label{mm}\left\{\begin{array}{l}
\left(A_i^*+\bar\lambda_i\right)A_{i,j}-A_{i,j}\left(A_j^*+\bar\lambda_j\right)
+\left(B_i^*+\bar\mu_i\right)B_{i,j}-B_{i,j}\left(B_j^*+\bar\mu_j\right)\\
 \ \ \ \ \ \ \ \ \ \ \ \ \ \ \ \ \ \ \ \ \ \ \ \ \ \ \ \ \ \ \ \ \ \ \ \ \ \ \ \ \ \ \ \ \ \ \ \ \ \ \ \ \ \ \ \ \ \ \ \ \ \ \ \ \ \ \ +Q_1(i,j)+Q_2(i,j)_2=x_ix_j^*\\
-\left(B_i+\mu_i\right)A_{i,j}+A_{i,j}\left(B_j+\mu_j\right)+
\left(A_i+\lambda_i\right)B_{i,j}-B_{i,j}\left(A_j+\lambda_j\right)+Q_3(i,j)=0\end{array}
\right.
\end{equation}
where $Q_1(i,j)$ (resp. $Q_2(i,j)$) is a quadratic expression depending on the $A_{\alpha,\beta}$'s 
(resp. in the $B_{\alpha,\beta}$)
and $Q_3(i,j)$ is bilinear in $A_{\alpha,\beta}$'s and $B_{\alpha,\beta}$.
For $\tau>0$ small enough, with the same arguments given in the preceding paragraph,
  the implicit function theorem implies 
\begin{lem}\label{approx}
The equations (\ref{mm}) has
a solution $A_{i,j},B_{i,j}\in \mcM_{n_i,n_j}(\C), i<j$ which depends smoothly on $\zeta\in \mcO_\Pg$,
moreover we have that
$$\sum_{i<j}\|A_{i,j}\|^2+\|B_{i,j}\|^2=O\left(\frac{1}{|\textbf{q}|^2}\right).$$
And the derivative of
the map $\zeta\mapsto (A_{i,j},B_{i,j})$ is bounded by $O\left(\frac{1}{|\textbf{q}|^2}\right).$
\end{lem}
Then we obtain $z_0(\zeta)=(A^0,B^0,x^0,0)\in \mcM_n(\C)\times\mcM_n(\C)\times\C^n\times\left(\C^n\right)^* $
an almost solution of the moment map equation :
$$\left\{\begin{array}{l}
\left[A^0,B^0\right]=0\\
2i\mu(z_0)-t=O\left(\frac{1}{|\textbf{q}|^2}\right)\end{array}
\right.$$
More precisely, the off block diagonal terms of the moment map equations are zero.
We will now used an argument that we learned in a paper of S. Donaldson \cite{Dk}[Proposition 17] :
we will find $h=ik$ a Hermitian matrix such that if
$$z_h=e^{ik}.z_0=(e^hA^0e^{-h},e^hB^0e^{-h},e^h.x^0,0)$$ then
$2i\mu(z_h)-t\Id=0$ and $\mu_\C(z_h)=0$ (this latter condition being obvious).

By the induction hypothesis $ii)$ and if $\tau$ is small enough and $R$ large enough then we have
$$\forall \eta\in \Ku_n,\ \|l_{z_0}(\eta)\|\geq C|\eta|$$
the constant $C$ being uniform on $\mcO_\Pg$.
Hence if $h=i\eta$  with $$\|k\|\le \delta:=\min\left\{1,\frac{Ce^{-2}}{4|z_0|}\right\}$$
then
$$\forall \eta\in \Ku_n,\ \|l_{z_h}(\eta)\|\ge \frac{C}{2}|\eta|$$
So as soon as we have
$\left|\mu(z_0)-\frac{t}{2i}\Id \right|< \frac{4}{C^2}\delta$, the proposition 17 in \cite{Dk}
 furnishes a $h=ik$ with
$\mu(e^h.z_0)=\frac{t}{2i}\Id$ with
$$\|h\|\le \frac{4}{C^2} \left|\mu(z_0)-\frac{t}{2i}\Id \right|.$$
But when $R$ is large enough, the condition 
$\left|\mu(z_0)-\frac{t}{2i}\Id \right|< \frac{4}{C^2}\delta$ is satisfied,
 hence there is $h=ik$ a Hermitian matrix such that $2i\mu(z_h)-t\Id=0$ and $\mu_\C(z_h)=0$.

We need to recall how $h$ is found.
For $z\in\M$ we have a linear map
$l_z:\Ku_n\rightarrow T_z\M\simeq \M$ and $l^*_z$ its adjoint; by definition of the moment map we have
$l_z^*=d\mu(z)\circ I$. The endomorphism $Q_z$ of $\Ku_n$ is given by $Q_z=l_z^*l_z$.
 Then for every $h=ik$, with $|k|<\delta$, 
$Q_{z_h}$ is invertible and $Q_{z_h}^{-1}$ has a operator norm bounded by $4C^{-2}$.
Let $a(z)=Q_z^{-1}\left(\mu(z)-\frac{t}{2i}\Id \right)$, we follow the maximal solution of the equation
\begin{equation}\label{EDO}
 \frac{dz}{ds}=-il_z(a(z)), 
\end{equation}
starting from $z_0$ at $s=0$.
 By definition we have
$$\frac{d\mu(z_s)}{ds}=-\left(\mu(z_s)-\frac{t}{2i}\Id \right)$$ hence
$$\mu(z_s)-\frac{t}{2i}\Id =e^{-s}\left(\mu(z_0)-\frac{t}{2i}\Id \right)\ ; $$
in fact $z_s=g_s.z_0$ where 
$$\frac{dg_s}{ds}=ia(z_s).g_s,\ g_s\in GL_n(\C)$$
The arguments of \cite{Dk} insures that the maximum solution of (\ref{EDO}) is defined on $[0,+\infty[$.
and if $g_s=e^{\eta_s}e^{h_s}$ where $\eta_s\in \Ku_n$ and $h_s$ is Hermitian,
then $|h_s|\le \delta$.
So that we also get :
 $$\|\dot g_s\|\le \frac{4}{C^2} \left|\mu(z_0)-\frac{t}{2i}\Id \right|e^{-s}e^\delta$$
hence $g_\infty=\lim_{s\to+\infty}g_s$ exists and
\begin{equation}\label{ginfini}
\|g_\infty-\Id\|\le \frac{4e^\delta}{C^2} \left|\mu(z_0)-\frac{t}{2i}\Id \right|
=O\left(\frac{1}{|\textbf{q}|^2}\right).
\end{equation}
We clearly have
$2i\mu(g_\infty.z_0)=t\Id$ and $h=ik$ is given by
$e^{2h}=g_\infty^*g_\infty$ i.e. the polar decomposition of 
$g_\infty$ is $g_\infty=e^{\eta_\infty}e^h$.
Moreover if $s\ge 0$, then the operator norm of
$l_{z_s}Q_{z_z}^{-1}$ remains less than $2/C$ 
hence
\begin{equation}\label{gz}
\|g_\infty.z_0-z_0\|\le \frac{2}{C} \left|\mu(z_0)-\frac{t}{2i}\Id \right|
=O\left(\frac{1}{|\textbf{q}|^2}\right).
\end{equation}
The implicit function theorem told us that $h$ depends smoothly    on $z_0$ hence on $\zeta\in \mcO_\Pg$, indeed
$$\left.\frac{d}{dt} \right|_{t=0}\mu(e^{tik}.z)=Q_{z}(k).$$
This map will be called : $\zeta\in \mcO_\Pg\rightarrow h(\zeta)\in i\Ku_n$.
The following lemma gives an estimate of the size of the derivative of $h$
\begin{lem}\label{dz}
Let $\textbf{v}\in T_\zeta\mcO_\Pg$  be a vector of unit length then
$$\|dh(\zeta).\textbf{v}\|=O\left(\frac{1}{|\textbf{q}|^2}\right).$$
\end{lem}
\proof 
Let $\dot z_0=dz_0(\zeta_0).\textbf{v}$ and $\dot h=dh(\zeta).\textbf{v}$, we also let $v\in \M$ be the
vector $v:=(\delta A,\delta B,\delta x,0)$ where if 
$$\textbf{v}=\left(\left((\delta\lambda_1,\delta\mu_1),...,(\delta\lambda_k,\delta\mu_k)\right),
\left(\delta A_1,...,\delta A_k\right), \left(\delta B_1,...,\delta B_k\right),(\delta x_1,..\delta x_k)\right)$$ then 
$\delta A$ (resp. $\delta B$) is a block diagonal matrix with diagonal 
$(\delta A_1+\delta\lambda_1\Id_{n_1},...,\delta A_k+\delta\lambda_k\Id_{n_k})$
(resp. $(\delta B_1+\delta\mu_1\Id_{n_1},...,\delta B_k+\delta\mu_k\Id_{n_k})$) and 
$\delta x=((\delta x_1,...,\delta x_k)$.

We have $$d\mu(z_h).\left(D\exp(h)\dot h.z_0
+e^h.\dot z_0\right)=0$$
Recall that : $$D\exp(h)\dot h=\frac{e^{\ad h}-\Id}{\ad h}.\dot h.e^h\ .$$
Let $i\dot\eta$ be the Hermitian part of $D\exp(h)\dot h$ and $\dot \xi$ be its skew Hermitian part.
Then
$$d\mu(z_h).(D\exp(h)\dot
h.z_0)=d\mu(z_h)(il_{z_h}\dot\eta)+d\mu(z_h)(l_{z_h}\dot\xi)=Q_{z_h}(\dot\eta).$$
Moreover from the construction of $z_0$ and the lemma (\ref{approx}), we obtain easily that
$$d\mu(z_0)(\dot z_0)=O\left(\frac{1}{|\textbf{q}|^2}\right) \ ;$$
 and $$\dot z_0=v+O\left(\frac{1}{|\textbf{q}|^2}\right)|v|.$$
 So if $k\in U(n)$ is such that $g_\infty=ke^h$ then
\begin{equation*}\begin{split}
\Ad\left(k\right)d\mu(z_h).\left(e^h.\dot z_0\right)&=d\mu(g_\infty.z_0).(g_\infty.\dot z_0)\\
&=d\mu(g_\infty.z_0).\left(\right(g_\infty-\Id\left).\dot z_0\right)
+d\mu(g_\infty.z_0-z_0).\dot z_0+d\mu(z_0)(\dot z_0) 
\end{split}
\end{equation*}
 Hence
 $$Q_{z_h}(\dot\eta)+d\mu(z_h).\left(k^{-1}.\right(g_\infty-\Id\left).\dot
 z_0\right)=O\left(\frac{1}{|\textbf{q}|^2}\right) \ .$$
 We now make the scalar product of this quantity with $\dot\eta$ and we obtain :
 \begin{equation*}\begin{split}
\|l_{z_h}(\dot \eta)\|^2&\leq O\left(\frac{1}{|\textbf{q}|^2}\right)\|\dot \eta\|-\langle
l_{z_h}(\dot\eta),k^{-1}.I\,\left(g_\infty-\Id\right).\dot z_0\rangle\\
&\leq O\left(\frac{1}{|\textbf{q}|^2}\right)\left( \|\dot \eta\|+\|l_{z_h}(\dot \eta)\|\right)
\end{split}\end{equation*}
 But our construction gives that
 $$\|\dot\eta\|\le \frac{2}{C} \|l_{z_h}(\dot \eta)\|\ ,$$
 hence we  obtain :
  $$\|\dot\eta\|\le \frac{2}{C} \|l_{z_h}(\dot \eta)\|=O\left(\frac{1}{|\textbf{q}|^2}\right).$$
Now $\dot h$ is a Hermitian matrix and $\|h\|=O\left(|\textbf{q}|^{-2}\right)$ hence by definition of $\dot\eta$ and $\dot\xi$, we have 
 $$ \dot h-i\dot\eta=O\left(|\textbf{q}|^{-2}\right).$$
Hence the lemma.
\endproof
We note that it is straightforward to verify the point $ii)$ at $z_h$ because by construction 
$$\forall \eta\in \Ku_n,\ l_{z_h}(\eta)|\geq \frac{C}{2}\|\eta\|.$$

We have build a map $f_\Pg$ from $\mcO_\Pg$ to $\LL_t$ whose value at
 a point $\zeta=\left(\textbf{q},\textbf{A},\textbf{B},\textbf{x}\right)\in \mcO_\Pg$ is the $z_h$
 constructed
 before. This map is $U(n_1)\times U(n_2)\times...\times U(n_k)$-equivariant hence it induces a
 map
 $$\psi_\Pg\,:\, \mcO_\Pg/(U(n_1)\times U(n_2)\times...\times U(n_k))\rightarrow \Hi{n}$$
 We remark that adjusting $\epsilon,R,\tau$, we have
 $$\mcC_\Pg\subset\mcO_\Pg/(U(n_1)\times U(n_2)\times...\times U(n_k))\subset \Hi{\Pg}\times
 \left(\C^2\right)^k_0.$$
 Where the last inclusion is an isometry if  $\Hi{\Pg}\times
 \left(\C^2\right)^k_0$ is endowed with the product metric.
 We now want to compare the metric $\psi_\Pg^*g_N$ and the product metric on 
 $\Hi{\Pg}\times
 \left(\C^2\right)^k_0$.
 Let $\textbf{v}$ be a vector of $T_\zeta\mcO_\Pg$ which is orthogonal to the $U(n_1)\times U(n_2)\times...\times U(n_k)$ 
 orbit of $\zeta$. As before, we defined $f(\zeta)=z_h=e^h .z_0$, $\dot h$, $v$..

 Recall that we have denote by $\Pi_{q}$ the orthogonal projection onto  $\im l_q$.
 Hence we need to compare
 $$\|\left(\Id-\Pi_{z_h }\right).df_\Pg(\zeta).\textbf{v}\|^2=\|df_\Pg(\zeta).\textbf{v}\|^2-\|\Pi_{z_h }.df_\Pg(\zeta).\textbf{v}\|^2\ \ {\rm and \ }\|\textbf{v}\|^2.$$
 But $$df_\Pg(\zeta).\textbf{v}=l_{z_h}(\dot\xi)+il_{z_h}(\dot\eta)+e^h.\dot z_0 $$
 Hence
 $$\left(\Id-\Pi_{z_h }\right).df_\Pg(\zeta).v=\left(\Id-\Pi_{z_h }\right).\left(il_{z_h}(\dot\eta)+e^h.\dot z_0\right).$$
 But we have already seen that
 $$\|l_{z_h}(\dot\eta)\|^2=O\left(\frac{1}{|\textbf{q}|^4}\right).$$
but $il_{z_h}(\dot\eta)$ is orthogonal to $T\LL_t$ hence to the range of $\Pi_{z_h}$, we also have 
\begin{equation*}\begin{split}
\langle \left(\Id-\Pi_{z_h }\right)\left(il_{z_h}(\dot\eta)\right),\left(\Id-\Pi_{z_h }\right).(e^h.\dot z_{0})\rangle
 &=\langle \left(\Id-\Pi_{z_h }\right).il_{z_h}(\dot\eta),(e^h.\dot z_{0})\rangle\\
&=\langle il_{z_h}(\dot\eta),(e^h.\dot z_{0})\rangle\\
&=-\left\langle \dot\eta,d\mu(z_h)\left(e^h.\dot z_0\right)\right\rangle\end{split}\end{equation*}
But
$$d\mu(z_h)\left((df_\Pg(\zeta).\textbf{v}\right)=0=d\mu(z_h)\left(il_{z_h}(\dot\eta)+e^h.\dot z_0\right)$$
so that 
\begin{equation*}\begin{aligned}\langle \left(\Id-\Pi_{z_h }\right)\left(il_{z_h}(\dot\eta)\right),\left(\Id-\Pi_{z_h }\right).(e^h.\dot z_{0})\rangle
 &=-\left\langle \dot\eta,d\mu(z_h)\left(e^h.\dot z_0\right)\right\rangle \\
&=\left\langle \dot\eta,d\mu(z_h)\left(il_{z_h}(\dot\eta)\right)\right\rangle=\left\| l_{z_h}(\dot\eta)\right\|\\
&=O\left(\frac{1}{|\textbf{q}|^4}\right).
\end{aligned}\end{equation*}
It remains now to estimate $$ \|\left(\Id-\Pi_{z_h }\right)(e^h.\dot z_{0})\|^2=
 \|e^h.\dot z_{0}\|^2-\|\Pi_{z_h }(e^h.\dot z_{0})\|^2.$$
 But $\Pi_{z_h}=l_{z_h}Q_{z_h}^{-1}l_{z_h}^*,$ and $$l_{z_h}^*\left(e^h.\dot z_{0}\right)=d\mu(z_h)(ie^h.\dot z_0).$$
We have already noticed that
$$I.\dot z_0=dz_0(\zeta).(I.\textbf{v})+w$$
 where $w=O(|\textbf{q}|^{-2}$.
So $$l_{z_h}^*\left(e^h.\dot z_{0}\right)=d\mu(z_h)\left(e^hdz_0(\zeta).(I\textbf{v})\right)+l_{z_h}^*(w)$$
And the proof of the lemma (\ref{dz}), furnishes a $w'=O(|\textbf{q}|^{-2})$ such that
$$d\mu(z_h)\left(e^hdz_0(\zeta).(I\textbf{v})\right)=d\mu(z_h)\left(w'\right)+O(|\textbf{q}|^{-2}$$
as the operator norm of $l_{z_h}Q_{z_h}^{-1}$ is bounded by $2/C$
we have obtained :
 $$\|\Pi_{z_h}(e^h.\dot z_0)\|^2=O\left(\frac{1}{|\textbf{q}|^4}\right).$$
 
 Hence we have obtain :
  \begin{equation*}\begin{split}
  \psi_\Pg^*g_N(\textbf{v},\textbf{v})&=\|e^h.\dot z_0\|^2+ O\left(\frac{1}{|\textbf{q}|^4}\right)\|v\|^2 \\
  &=\|\dot z_0\|^2+2\langle \dot z_0,h\dot z_0\rangle+O\left(\frac{1}{|\textbf{q}|^4}\right)\|v\|^2. \\
  \end{split}\end{equation*}
  By construction 
  $$|\dot z_0|^2=\|v\|^2+O\left(\frac{1}{|\textbf{q}|^4}\right)\|v\|^2$$
  And if
  $\textbf{v}=(\delta q,\delta A_1,\delta A_2,...,\delta A_k,\delta B_1,\delta B_2,...,\delta B_k,0)$
  and if $h_{i,j}$ are the block of $h$ of size $n_i\times n_j$ then
  \begin{equation*}\begin{split}\langle \dot z_0,h\dot z_0\rangle&= \langle v ,hv\rangle+O\left(\frac{1}{|\textbf{q}|^4}\right)\|v\|^2\\
  &=\sum_j \langle \delta A_j, \left[h_{j,j},\delta A_j\right]\rangle
  +\langle \delta B_j, \left[h_{j,j},\delta B_j\right]\rangle
  +O\left(\frac{1}{|\textbf{q}|^4}\right)\|v\|^2\\
 &=\sum_j \langle \left[\delta A_j, \delta A_j^*\right]+
  \left[\delta B_j, \delta B_j^*\right],h_{j,j}\rangle
  +O\left(\frac{1}{|\textbf{q}|^4}\right)\|v\|^2\\
  &= O\left(\frac{1}{\sigma^2|\textbf{q}|^2}\right)\|v\|^2
  \end{split}\end{equation*}
  according to the hypothesis $iii)$.
  
  In order to finish the proof we need to check the property $iii)$ at the point $z_h$.
 With what has been proved in the preceding paragraph, we only need to check that
  if $(\delta A,\delta B,0,0)$ is a unitary vector in the
   tangent space of $\LL_t$ at $z_h$ and orthogonal to $U(n)$ orbit
  of $z_h$ then $\left[\delta A_j, \delta A_j^*\right]+\left[\delta B_j, \delta B_j^*\right]$ is bounded.
  This is evident.
\section{Conclusion}
With the previous asymptotic of Nakajima's metric, we'll show that Nakajima's metric coincides with Joyce's
one ;  a way for proving such a result would be to verify the estimate \ref{QALE1}) for the orders $l\ge 1$ ;
this is probably possible with some extra work, however we'll give here a different proof which follows the
classical proof of the unicity for the solution of the Calabi-Yau problem. Moreover our argument gives a new
analytical result on mapping property of the Laplace operator on QALE space. For new results which extended
Joyce's ones and which go further than our result, there is a forthcoming work of A. Degeratu and R. Mazzeo \cite{DM}.

We have already seen that Kronheimer's classification of hyperk\"ahler ALE $4-$dimensional manifold implies
that on $\Hi{2}\simeq T^*\PP^1(\C)$ , Nakajima's metric is the Eguchi-Hansen metric. We are going to prove
our result by induction on $n$. Hence we now assume that up to a scaled factor, Joyce's and Nakajima's metrics
coincide on $\Hi{l}$ for all integer $l<n$.
We consider $g$ the Joyce's metric on $\Hi{n}$ and $\omega$ the K\"ahler form associated to $g$ (for the complex structure $I$) and for simplicity of forthcoming notation, we denote by $g'$
Nakajima's metric on $\Hi{n}$ with associated K\"ahler form $\omega'$.
\subsection{Comparison of the two metrics} 
The second group of cohomology of $\Hi{n}$ has dimension $1$ and a cycle  dual to a basis of $H^2(\Hi{n},\R)$
is given by the image of a holomorphic map
$f_n\,:\, \PP^1(\C)\rightarrow \Hi{n}$ such that if $\pi\,:\, \Hi{n}\rightarrow \left(\C^2\right)^n_0/S_n$
then the image of $f_n$ is $\PP^{1}(\C)\simeq\pi^{-1}\{((0,0),q')\}S_n$ for 
$q'\in\left(\C^2\setminus\{0\}\right)^{n-2}_0$. We can assume that 
$$\int_{\PP^1(\C)} f_n^*\omega=\int_{\PP^1(\C)}  f_n^*\omega'$$ 
Moreover, for each partition $\Pg$ of $\{1,2,...,n\}$, we have on $\mcC_\Pg\subset \Hi{\Pg}\times V_\Pg$
$$\psi^*_\Pg g=g_\Pg+\eucl+O\left(\frac{1}{\sigma^2\rho^2}\right)$$ and
$$\psi^*_\Pg g'=g'_\Pg+\eucl+O\left(\frac{1}{\sigma^2\rho^2}\right)$$
where $g_\Pg$ (resp. $g_\Pg'$) is the sum of the Joyce's (resp. Nakajima's) metric on $\Hi{\Pg}\simeq \Hi{n_1}\times
\Hi{n_2}\times ...\times \Hi{n_k}$.
However $f_n(\PP^1(\C))$ is homologous to 
$C_{\Pg,j,v}
:=\psi_{\Pg}\{(y_1,...,y_{n_j -1}\}\times f_{n_j}(\PP^1(\C))\times\{y_{n_j+1},..,y_{n_k}\}\times \{v\}$
where $y_j\in \Hi{n_j}$ and $v\in V_\Pg$, let $g_j$ (resp. $(g_j'$) be the Joyce's (resp. Nakajima's) metric on  $\Hi{n_j}$
and $\omega_j$ (resp. $\omega_j'$) its K\"ahler form ; that is to say $g_\Pg=g_1+g_2+...++g_k$ and
 $g'_\Pg=g'_1+g'_2+...+g'_k$. We have
 
\begin{equation*}
\begin{split}
\int_{C_{\Pg,j,v}}\omega&=\int_{f_{n_j}(\PP^1(\C))}\omega_j+O\left(\frac{1}{\sigma^2\rho^2}\right)\\
&=\int_{C_{\Pg,j,v}}\omega'\\
&=\int_{f_{n_j}(\PP^1(\C))}\omega'_j+O\left(\frac{1}{\sigma^2\rho^2}\right)
\end{split}\end{equation*}
In particular letting $\|v\|$ going to $\infty$, we obtain
$$\int_{f_{n_j}(\PP^1(\C))}\omega'_j=\int_{f_{n_j}(\PP^1(\C))}\omega_j$$
Our induction hypothesis yields that $g_j=g_j'$ for all $j$, and
eventually, we have proved that
$$g-g'=O\left(\frac{1}{\sigma^2\rho^2}\right)$$
\subsection{Coincidence of Joyce's and Nakajima's metrics}
Following the classical proof of the unicity of the solution of the Calabi-Yau problem, we would like to find
a good function $\phi$ such that
$$\omega-\omega'=i\partial\bar\partial\phi.$$
However it is not easy because the weight $\sigma^{-2}\rho^{-2}$ is critical in Joyce's analysis on QALE
manifold. To circumvent this difficulty, we remark that both metrics $g$ and $g'$ have a $\bS^1$ invariance
property coming from the diagonal action of $\bS^1$ on $\left(\C^2\right)^n_0/S_n$. For Joyce's metric it
comes from the unicity result of the QALE K\"ahler Einstein metric asymptotic to $\left(\C^2\right)^n_0/S_n$.
For Nakajima's metric, the action of $\bS^1$ on $\M$ is the following :
if $e^{i\theta}\in \bS^1$ and if $z=(A,B,x,0)\in \LL_t$ then
$e^{i\theta}.z:=\left(e^{i\theta}A,e^{i\theta}B,e^{i\theta}x,0\right)\in \LL_t$. And this action is isometric.
This $\bS^1$ action is holomorphic for the complex structure $I$ but not for the complex structures $J$ and 
$K$. Let $X$ be the $g$ or $g'$ Killing field associated to the infinitesimal action of  
$\eta= i/2$.
Then $X$ has linear growth on $\Hi{n}$ that is to say
there is a constant $c$ such that 
$$X(z)\leq c(\rho(z)+1).$$
Moreover if $\omega_1$ is the K\"ahler form of $(g, J)$, $\omega_2$ is the K\"ahler form of $(g, K)$ and 
$\omega_1'$ and $\omega'_2$ are the corresponding form associated to the metric $g'$ then
$$\omega_1=d(i_X\omega_2) {\rm \ and\ }\omega_1'=d(i_X\omega'_2)$$
Hence if we let $$\beta=i_X\omega_2-i_X\omega'_2$$ then we have
$$\omega_1-\omega_1'=
d\beta {\rm \ and\ }\beta=O\left(\frac{1}{\sigma^2\rho}\right)$$

We work now in the K\"ahler manifold $(\Hi{n},g, J)$, the following analytical result is the key point of our
proof :
\begin{prop}\label{ana}
There is a $(0,1)$-form $\alpha$ on $\Hi{n}$ such that
$$\alpha=O\left(\frac{\log(\rho+2)}{\rho}\right)$$
and 
$$\beta^{0,1}=\Delta_{\bar\partial}\alpha=\bar\partial\bar\partial^*\alpha+
\bar\partial^*\bar\partial\alpha.$$
\end{prop}
We first explain how we can prove that $\omega_1=\omega_1'$ with this proposition. This proposition will be
proved in the next subsection.

The $1$-form $\Phi=\bar\partial^*\bar\partial\alpha$ satisfies $\bar\partial\beta^{0,1}=0=\bar\partial\Phi$
and
$\bar\partial^*\Phi=0$.
Moreover the metric $g$ has by definition bounded geometry, hence we have the 
following uniform in $x\in \Hi{n}$ local elliptic estimate :
\begin{equation*}
\begin{split}
\|\Phi\|_{L^2(B(x,1))}=\|\bar\partial^*\bar\partial\alpha\|_{L^2(B(x,1))}&\leq c
\|\Delta_{\bar\partial}\alpha\|_{L^2(B(x,2))}+c'\|\alpha\|_{L^2(B(x,2))}\\
&\leq O\left(\frac{\log(\rho+2)}{\rho}\right)\end{split}\end{equation*}
But $\Phi$ being harmonic we also have a uniform estimate 
$$|\Phi(x)|\leq c\|\Phi\|_{L^2(B(x,1)}.$$
Hence we obtain that 
$$\Phi=O\left(\frac{\log(\rho+2)}{\rho}\right).$$
But the Ricci curvature of $g$ is zero hence the Bochner formula 
and the Kato inequality implies that $|\Phi|$ is a subharmonic function hence $\Phi$ is zero by the maximum
principle. And we get $\beta^{0,1}=\bar\partial\bar\partial^*\alpha$, the same argument shows that we can
find a $(1,0)$-form $\tilde\alpha$ such that $\beta^{1,0}=\partial\partial^*\tilde\alpha$. Hence if we let
$$i\phi=\bar\partial^*\alpha-\partial^*\tilde\alpha$$ then we have
$$d\beta=i\partial\bar\partial\phi.$$
Again the same argument as before using the fact that $g$ has bounded geometry, implies that
$$\phi=O\left(\frac{\log(\rho+2)}{\rho}\right).$$

Both $\omega_1$ and $\omega_1'$ are K\"ahler Einstein with zero scalar curvature hence there is a pluriharmonic
function $f$ such that
$$\omega_1^m=e^f\left(\omega_1'\right)^m$$
But we also have $$f=O\left(\frac{1}{\sigma^2\rho^2}\right).$$
By the maximum principle we deduce that $f=0$.
We finish the proof with a classical argument : the function $\phi$ is subharmonic for the metric $g$
\cite{besse}[expos\'e VI, lemma 1.6] and decay at infinity hence by the maximum principle 
$\phi$ is negative ; but reversing the role of $g$ and $g'$,
 $-\phi$ is also subharmonic for the metric $g'$ and $-\phi$ is positive and decay at infinity hence $\phi$ is zero.

\subsection{Proof of the analytical result}We first remark that because $(\Hi{n},g)$ is asymptotic to the
Euclidean metric on $\left(\C^2\right)^n_0/S_n$, we have
$$\lim_{r\to\infty} \frac{\vol B(x,r)}{r^d}=\frac{w_d}{n!}$$
where $d=4(n-1)$ is the real dimension of $\Hi{n}$ and 
$w_d$ is the volume of the unit ball in $\R^d$.
The Bishop-Gromov inequality tolds us that for any point $x\in \Hi{n}$ :
$$\frac{w_d r^d}{n!}\le \vol B(x,r)\leq w_d r^d.$$

The result of P. Li and S-T. Yau implies that the Green kernel $G$ of the metric $g$ (that is to say the Schwartz kernel
of the operator $\Delta^{-1}$) satisfies \cite{LiYau} :
$$G(x,y)\le \frac{c}{d(x,y)^{d-2}}.$$

Moreover because $g$ is Ricci flat, the Hodge-deRham operator acting on $1$ forms is the rough Laplacian :
$$\forall v\in C^\infty_0(T^*\Hi{n}),\ \ \Delta=dd^*+d^*d=\nabla^*\nabla$$
Hence the Kato inequality implies that if $\vec G(x,y)$ is the Schwartz kernel of the operator
$\Delta^{-1}$, then it  satisfies
$$|\vec G(x,y)|\leq  G(x,y)\le \frac{c}{d(x,y)^{d-2}}.$$

The proposition (\ref{ana}) will be a consequence of the following lemma
 \begin{lem} If $f\in L^\infty_{\rm loc}(\Hi{n})$  is a non negative function which satisfies 
 $$f=O\left(\frac{1}{\sigma^2\rho}\right),$$
 then $$u(x)=\int_{\Hi{n}} \frac{f(y)}{d(x,y)^{d-2}}dy$$
 is well defined and satisfies
 $$u=O\left(\frac{\log(\rho+2)}{\rho}\right).$$
 \end{lem}
 \proof Let $o\in \Hi{n}$ be a fixed point. And we can assume that $\rho(x)=d(o,x)$. 
 We remark that $u$
  is well defined indeed there is a constant $c$ such that for $R>1$ then
  $$\int_{B(o,R)} f\leq c R^{d-3}.$$
  As a matter of fact, the function $\frac{1}{\sigma^2\rho}$ is asymptotic to a homogeneous function
  $h$ of degree $-3$ on $\left(\C^2\right)^n_0/S_n$
  $h(r\theta)=r^{-3}\bar h(\theta)$ where 
  $\bar h$ is a positive function on $\bS^{2d-1}/S_n$ ; this function $\bar h$ is singular
  on the singular  locus of $\bS^{2d-1}/S_n$; we call $\Sigma$ this singular locus.
   but $\bar h$ behaves like $d(.,\Sigma)^{-2}$ near $\Sigma$ but the real co dimension of $\Sigma$ is 
   $4$ hence $\bar h$ is integrable on $\bS^{2d-1}/S_n$ and
   we have
   $$\lim_{R\to \infty} R^{3-d}\int_{B(o,R)}\frac{1}{\rho\sigma^2}= \frac{1}{d-3}\int_{\bS^{2d-1}/S_n} \bar h.$$
   
In order to finish our estimate, we must find a constant $c$ such that if $\rho(x)\ge 10$ then 
$$F(x)=\int_{\Hi{n}}\frac{1}{d(x,y)^{d-2}}\frac{1}{\rho(y)\sigma(y)^2} dy\le c
   \ \frac{\log\rho(x)}{\rho(x)}$$
We decompose

\begin{equation}\label{decomp}\Hi{n}=\left(B(o,2\rho(x))\setminus B(x,\rho(x)/2)\right)\cup B(x,\rho(x)/2)\cup \left(\Hi{n}\setminus
B(o,2\rho(x))\right)\ ,\end{equation}
then we have
$F=F_1+F_2+F_3$ where $F_i$ is the integral of $d(x,y)^{2-d}\rho^{-1}\sigma^{-2}$ on the $i^{\rm th}$ region
of the decomposition (\ref{decomp}).
The first and the last integrals are easy to estimate :
$$F_1(x)\le \left(\frac{2}{\rho(x)}\right)^{d-2}\int_{B(o,2\rho(x))}\frac{1}{\rho(x)\sigma(x)^2}\le C
\frac{1}{\rho(x)}.$$
Concerning $F_3$ we have
\begin{equation*}
\begin{split}
F_3(x)&=\int_{\Hi{n}\setminus B(o,2\rho(x))}\frac{2^{d-2}}{d(x,y)^{d-2}}\frac{1}{\rho(y)\sigma(y)^2} dy\\
&\leq \int_{\Hi{n}\setminus B(o,2\rho(x))}\frac{2^{d-2}}{\rho(y)^{d-1}}\frac{1}{\sigma(y)^2} dy\\
&\leq \sum_{k=1}^\infty \int_{B(o,2^{k+1}\rho(x))\setminus B(o,2^k\rho(x))}
\frac{2^{d-2}}{\rho(y)^{d-1}}\frac{1}{\sigma(y)^2} dy\\
&\leq \sum_{k=1}^\infty \frac{1}{\left(2^k\rho(x)\right)^{d-2}}
\int_{B(o,2^{k+1}\rho(x))}\frac{1}{\rho(y)\sigma(y)^2} dy\\
&\leq C\sum_{k=1}^\infty \frac{1}{\left(2^k\rho(x)\right)^{d-2}} \left(2^{k+1}\rho(x)\right)^{d-3}\\
&\leq C'\frac{1}{\rho(x)}.
\end{split}
\end{equation*}
It remains to estimate $F_2$ :
We have
$$F_2(x)\leq \frac{2}{\rho(x)} \int_{B(x,\rho(x)/2)} \frac{1}{d(x,y)^{d-2}}\frac{1}{\sigma(y)^2}dy.$$
Let $V(\tau)=\int_{B(x,\tau)}\frac{1}{\sigma(y)^2}dy$ and note $dV$ the Riemann-Stieljes measure associated
to the increasing function $V$. We have
\begin{equation}\label{SJ}
\begin{split}
\int_{B(x,\rho(x)/2)} \frac{1}{d(x,y)^{d-2}}\frac{1}{\sigma(y)^2}dy&=
\int_{0}^{\rho(x)/2}\frac{1}{\tau^{d-2}}dV(\tau)\\
&=\frac{V(\rho(x)/2)}{\left(\rho(x)/2\right)^{d-2}}+(d-2)\int_0^{\rho(x)/2}\frac{V(\tau)}{\tau^{d-1}}d\tau.
\end{split}
\end{equation}
We will estimate $V$ :
if we note $S$ the pull back to $\Hi{n}$ of the singular locus of $\left(\C^2\right)^n_0/S_n$ 
and $\mcO=\{y\in \Hi{n},{\rm\ such\ that\ } \sigma(y)\ge 2\}$
then we 
have $V(\tau)=V_1(\tau)+V_2(\tau)$
where $V_1$ is the integral over $B(x,\tau)\cap\mcO$ and $V_2$ is the integral over $B(x,\tau)\setminus\mcO$.

$V_1$ is easy to estimate because on this region, $\sigma^{-2}$ is bounded hence
\begin{equation}\label{pres}
 V_1(\tau)\le C\vol \left(B(x,\tau)\cap \mcO\right) \leq C\min\{\tau^d,\tau^{d-4}\}.
\end{equation}

Outside $\mcO$ the metric is quasi-isometric to the Euclidean metric and we can estimate
$V_2$ by a similar integral on $\left(\C^2\right)^n_0/S_n$.
Let $$D=\{q\in\left(\C^2\right)^n_0, {\rm\ such\ that\ } \forall i\not=j\ \ |q_i-q_j|\geq |q_1-q_2|\}$$
and $D'=\{q\in D, |q_1-q_2|\geq 1\}$.
$D$ is a fundamental domain for the action of $S_n$ on $\left(\C^2\right)^n_0$ and if $\bar x\in D$ is such
that $S_n.\bar x=\pi (x)$ then
$$V_2( \tau)\leq\frac{C}{n!}\sum_{\gamma\in S_n}\int_{D'\cap B(\gamma\bar x,\tau)}\frac{1}{|q_1-q_2|^2}dq.$$
We give three different estimates for $V_2$ according to the relative size of $\sigma(x)$ and $\tau$:
\begin{enumerate}
\item If $\sigma(x)\leq 3/2$ then for $\tau\in[0,1/2]$ we have $V_2(\tau)=0$.
\item if $\sigma(x)\geq 3/2$ then for $\tau\leq \sigma(x)/2$ then 
$$V_2(\tau)\leq \frac{C}{\sigma(x)^2}\tau^d$$
\item and finally if  $\tau\geq \sigma(x)/2$ then
there is a point $z\in S$ such that $d(x,z)=\sigma(x)-1$ and if $\bar z\in D$ such that
$S_n\bar z=\pi(z)$
then
$$\int_{D'\cap B(\gamma\bar x,\tau)}\frac{1}{|q_1-q_2|^2}dq
\leq \int_{D'\cap B(\gamma\bar z,3\tau)}\frac{1}{|q_1-q_2|^2}dq\leq C\tau^{d-2}.$$
\end{enumerate}
Now, with the estimate (\ref{pres}), it is easy to show that in (\ref{SJ})
 the part coming from $V_1$ is bounded ; 
  concerning the part coming from $V_2$, when $\sigma(x)\le 3/2$, we get
\begin{equation*}
\int_{0}^{\rho(x)/2}\frac{1}{\tau^{d-2}}dV_2(\tau)
\leq C+(d-2)\int_{3/2}^{\rho(x)/2} \frac{C\tau^{d-2}}{\tau^{d-1}}d\tau=C'+C\log\rho(x)\ ,
\end{equation*}
and when $\sigma(x)\ge 3/2$, we obtain 
\begin{equation*}\begin{split}
\int_{0}^{\rho(x)/2}\frac{1}{\tau^{d-2}}dV_2(\tau)&\leq C+(d-2)\int_0^{\sigma(x)/2} \frac{C\tau^{d}}{\tau^{d-1}\sigma(x)^2}d\tau+(d-2)\int_{\sigma(x)/2}^{\rho(x)/2} \frac{C\tau^{d-2}}{\tau^{d-1}}d\tau\\
&=C'+C\log\left(\frac{\rho(x)}{\sigma(x)}\right)\ 
\end{split}\end{equation*}
Hence the result.
 \endproof


\end{document}